\documentclass[12pt,a4paper]{article}
\usepackage[utf8x]{inputenc}
\usepackage[english]{babel}
\usepackage[T1]{fontenc}
\usepackage{amsmath}
\usepackage{amsfonts}
\usepackage{amssymb}
\usepackage{makeidx}
\usepackage{graphicx}
\usepackage[left=3cm,right=2cm,top=3cm,bottom=2cm]{geometry}
\usepackage{setspace}
\usepackage{amsthm}
\usepackage{cite}

\usepackage[colorlinks=true,urlcolor=blue,
citecolor=red,linkcolor=blue,linktocpage,pdfpagelabels,
bookmarksnumbered,bookmarksopen]{hyperref}

\addto\captionsenglish{}

\newtheorem{definition}{Definition}[section]

\newtheorem{proposition}{Proposition}[section]

\newtheorem{theorem}{Theorem}[section]

\newtheorem{example}{Example}[section]

\newtheorem{lemma}{Lemma}[section]

\newtheorem{observation}{Remark}[section]

\newtheorem{corollary}{Corolary}[subsection]

\numberwithin{equation}{section}

\newcommand{\bo}{\begin{observation}}
\newcommand{\eo}{\end{observation}}
\newcommand{\bd}{\begin{definition}}
\newcommand{\ed}{\end{definition}}
\newcommand{\bp}{\begin{proposition}}
\newcommand{\ep}{\end{proposition}}
\newcommand{\bt}{\begin{theorem}}
\newcommand{\et}{\end{theorem}}
\newcommand{\bc}{\begin{corollary}}
\newcommand{\ec}{\end{corollary}}
\newcommand{\bl}{\begin{lemma}}
\newcommand{\el}{\end{lemma}}
\newcommand{\be}{\begin{example}}
\newcommand{\ee}{\end{example}}
\newcommand{\beq}{\begin{equation}}
\newcommand{\eeq}{\end{equation}}
\newcommand{\beqa}{\begin{equation*}}
\newcommand{\eeqa}{\end{equation*}}
\newcommand{\R}{\mathbb{R}}
\newcommand{\RN}{\mathbb{R}^{N}}
\newcommand{\N}{\mathbb{N}}
\newcommand{\Rdois}{\mathbb{R}^{2}}

\newcommand{\Ldois}{ L^{2}(\mathbb{R}^2) }
\newcommand{\Loito}{L^{\frac{8}{3}}(\Rdois)}

\newcommand{\Hum}{H^{1}(\mathbb{R}^2)}
\newcommand{\Ls}{ L^{s}(\mathbb{R}^2)}

\newcommand{\intR}{\displaystyle\int\limits_{\mathbb{R}^2}}
\newcommand{\un}{u_{n}}

\newcommand{\vn}{v_{n}}

\newcommand{\sn}{s_{n}}
\newcommand{\yn}{y_{n}}
\newcommand{\wn}{w_{n}}
\newcommand{\util}{\tilde{u}}
\newcommand{\until}{\tilde{u}_n}
\newcommand{\RA}{\rightarrow}

\newcommand{\ds}{\displaystyle\int\limits}
\newcommand{\intlog}{ \displaystyle\int\limits_{\mathbb{R}^2} \displaystyle\int\limits_{\mathbb{R}^2} \ln (|x-y|)u^2(x)u^2(y) dx dy}
\newcommand{\intLog}{ \displaystyle\int\limits_{\mathbb{R}^2} \displaystyle\int\limits_{\mathbb{R}^2} \ln (|x-y|)u^2(x)u(y)\varphi(y) dx dy}

\newcommand{\cqd}{\hfill $\rule{2mm}{2mm}$}
\usepackage{times, color, xcolor}
\addto\captionsportuguese{
}

\begin{document}

 	\title{The  Choquard logarithmic 
equation involving  a  nonlinearity   with  exponential growth  
		\thanks{The first author was supported    by  Coordination of Superior Level Staff Improvement-(CAPES) -Finance Code 001 and  S\~ao Paulo Research Foundation- (FAPESP), grant $\sharp $ 2019/22531-4,
while the second  author was supported by  National Council for Scientific and Technological Development -(CNPq),   grant $\sharp $ 307061/2018-3 and FAPESP  grant $\sharp $ 2019/24901-3.
			}}
	\author{
		Eduardo  de S. Böer \thanks{ E-mail address: eduardoboer04@gmail.com Tel. +55.51.993673377}  and Ol\'{\i}mpio H. Miyagaki \thanks{Corresponding author} \footnote{ E-mail address: ohmiyagaki@gmail.com, Tel.: +55.16.33519178 (UFSCar).}\\
		{\footnotesize Department of Mathematics, Federal University of S\~ao Carlos,}\\
		{\footnotesize 13565-905 S\~ao Carlos, SP - Brazil}\\ }
\noindent
				
	\maketitle

\noindent \textbf{Abstract:} In the present work we are concerned with the Choquard Logarithmic equation $-\Delta u + au + \lambda (\ln|\cdot|\ast |u|^{2})u = f(u) \textrm{ \ in \ } \mathbb{R}^2$, for $ a>0 $, $ \lambda >0 $ and a nonlinearity $f$ with exponential critical growth. We prove the existence of a nontrivial solution at the mountain pass level and a nontrivial ground state solution. Also, we provide these results under a symmetric setting, taking into account subgroups of $ O(2) $.

\vspace{0.5 cm}


\section{Introduction}

In this paper we are interested in studying standing wave solutions for the following Schrödinger-Poisson System,
\beq\label{i5}
\left\{ \begin{array}{rclcl}
i\psi_t -\Delta\psi + \tilde{V}(x)\psi + \gamma \omega \psi &= &0 \  &  \textrm{ in \ } &  \RN\times \R \\
\Delta \omega &=& |\psi|^2 & \textrm{ in \ } & \RN,
\end{array} \right. 
\eeq
where $ \psi : \RN \times \R \RA \mathbb{C} $ is the time-dependent wave function, $ \tilde{V}: \RN \RA \R $ is a real external potential and $ \gamma > 0 $ is a parameter. The function $ \omega $ represents an internal potential for a nonlocal self-interaction of the wave function $ \psi $. The usual ansatz $ \psi(x, t)=e^{-i\theta t}u(x) $, with $ \theta \in \R $, for standing wave solutions of (\ref{i5}) leads to
\beq\label{i1}
\left\{ \begin{array}{rclcl}
-\Delta u + V(x)u + \gamma \omega u & = &0  & \textrm{ in \ }& \RN \\
\Delta \omega &= &u^2 &\textrm{ in \ }& \RN,
\end{array} \right. 
\eeq
with $ V(x)= \tilde{V}(x) + \theta $. The second equation of (\ref{i1}) shows us that $ \omega : \RN \RA \R $ is determined only up to harmonic functions. In this point of view, it is natural to choose $ \omega $ as the Newton potential of $ u^2 $, i.e., $ \Gamma_N \ast u^2 $, where $ \Gamma_N $ is the well-known fundamental solution of the Laplacian
$$
\Gamma_N (x) = \left\{ \begin{array}{rcl}
\dfrac{1}{N(2-N)\sigma_N}|x|^{2-N} &\textrm{ if \ } &N\geq 3 ,\\
\dfrac{1}{2\pi}\ln |x|& \textrm{ if \ }& N= 2 ,
\end{array} \right.
$$
where $ \sigma_N $ denotes the volume of the unit ball in $ \RN $. With this formal inversion of the second equation in (\ref{i1}), as  it is detailed in \cite{[4]}, we obtain the following integro-differential equation
\beq\label{i2} 
-\Delta u + V(x) u + \gamma ( \Gamma_N \ast |u|^2) u = b|u|^{p-2}u,\ p> 2, \ b>0, \textrm{ \  \ in \ } \RN.
\eeq
The case $ N=3 $ has been extensively studied, due to its relevance in physics. Although this equation is called ``Choquard equation'', it has first studied by Fröhlich and Pekar in \cite{[12] , [11] , [22]}, to describe the quantum mechanics of a polaron at rest, in the particular case, when $ V(x)\equiv a > 0 $ and $ \gamma > 0 $. Then, it was introduced by Choquard in 1976, to study  an electron trapped in its hole. The local nonlinear terms on the right side of equation (\ref{i2}), such as $ b|u|^{p-2}u $, for $ b\in \R $ and $ p>2 $, usually appears in Schrödinger equations to model the interaction among particles. Since then, several variations of Choquard equations have been studied to model a series of phenomena.

Summarizing the discussion about the three dimensional case of equation (\ref{i2}), we refeer to some additional papers, and the references therein, \cite{[14] , [16] , [1] , [2] , [3] , [20] , [18]}. 

Once we turn our attention to the case $ N=2 $, we immediately see that the literature is scantier. In this case, we can cite the recent works of \cite{[6], [cjj], [10]}. In \cite{[6]}, the authors have proved the existence of infinitely many geometrically distinct solutions and a ground state solution, considering $ V: \Rdois \RA (0, \infty) $ continuous and $ \mathbb{Z}^2 $-periodic, $ \lambda > 0 $ and a particular case $ f(u)=b|u|^{p-2}u $, with $ b\geq 0 $ and $ p\geq 4 $. Here  because of the  periodic setting, the global Palais–Smale condition can fail, since  the corresponding functional become invariant under $ \mathbb{Z}^2 $-translations. Then, intending to fill the gap, Du and Weth \cite{[10]} studied equation (\ref{i2}) in the case $ V(x)\equiv a > 0 $, $ \lambda >0 $ and $ f(u)=|u|^{p-2}u $, with $ 2<p<4 $. They have proved the existence of a mountain pass solution and a ground state solution. Also, they verified that, if $ p\geq 3 $, both levels are equal and provided a characterization for them.  Finally, in \cite{[cjj]}, the authors dealt with the existence of stationary waves with prescribed norm considering $ \lambda \in \R $.


In the present paper we foccus on planar case 
\begin{equation} \label{P}
-\Delta u + au + \lambda (\ln|\cdot|\ast |u|^{2})u = f(u) \textrm{ \ in \ } \mathbb{R}^2,
\end{equation}
where $ a>0 $, $ \lambda >0 $ and $f: \R \RA [0, \infty) $ is continuous, with primitive $ F(s)=\int\limits_{0}^{s}f(t)dt $. 
We recall that Stubbe, in\cite{[21]},  set up a variational framework for (\ref{P}), with $f=0,$ within a subspace of $ \Hum $, where the associated functional is well-defined.

As one can note, the results contained in \cite{[6] , [10]} use the geometrical form of the nonlinearity, $ f(u)=b|u|^{p-2}u, $ to obtain the results without a bigger concern on the boundedness of the sequences involved. Our paper represents an extension of the mentioned results, since it takes into consideration a more general growth condition over $ f $, to be precise, we will consider an exponential critical growth to the nonlinearity $f$, introduced in the seminal paper \cite{[17]}. As we will point out bellow, it rises new difficulties while dealing with equation (\ref{P}).  To the best of our knowledge, it is the first paper combining logarithmic Choquard equations with exponential critical growth. 

We say that a funcion $ h $ has \textit{subcritical} exponential growth at $ +\infty $, if
$$
\lim\limits_{s\RA + \infty}\dfrac{h(s)}{e^{\alpha s^2}-1} = 0 \textrm{ \ , for all \ } \alpha >0 ,
$$
and we say that $ h $ has $ \alpha_0 $-\textit{critical} exponential growth at $ +\infty $, if
$$
\lim\limits_{s\RA + \infty}\dfrac{h(s)}{e^{\alpha s^2}-1} = \left\{ \begin{array}{ll}
0, \ \ \ \forall \ \alpha > \alpha_{0} \\
+\infty , \ \ \ \forall \ \alpha < \alpha_0
\end{array} \right. .
$$
As usual conditions while dealing with this kind of growth, found in works such as \cite{[19] , [9]}, we assume that $ f $ satisfies
 $$\left\{ \begin{array}{c}f: \R \RA [0, \infty) \ \mbox{is continuous and has critical exponential growth with } \ \alpha_0 = 4 \pi,   \\
 |f(s)|\leq C e^{4\pi s^2}, \forall \ s \in \R
 \end{array}\right. \leqno{(f_1)}$$
$$ \lim\limits_{s\RA 0} \dfrac{f(s)}{s}=0.  \leqno{(f_2)}$$
From $ (f_1) $ and $(f_2)$, given $ \varepsilon >0 $, $\tau>1,$ fixed, for all  $ p>2 $, we can find two constants $ b_1, b_2 > 0$ such that 
\begin{equation}\label{eq1}
f(s)\leq \varepsilon |s| +b_1 |s|^{p-1}(e^{\tau 4 \pi s^2}-1) \ , \ \ \ \forall \ s \in \R,
\end{equation} 
and
\begin{equation}\label{eq2}
F(s) \leq \dfrac{\varepsilon}{2}|s|^2 + b_2 |s|^p(e^{\tau 4\pi s^2}-1) \ , \ \ \ \forall s \in \R .
\end{equation}

As one can see, our new concern will be about how to guarantee the boundedness of an integral of the form $ \ds_{\Rdois} (e^{\alpha u_{n}^{2}}-1) dx $, with $ \alpha > 0 $, for a specific Cerami sequence $(\un)$. Later, we will present the condition that will make it possible. 

Now, associated with our problem we have the following functional $ I: \Hum \RA \R \cup \{ \infty\} $ given by
\begin{equation}\label{I}
I(u) = \dfrac{1}{2}\intR |\nabla u |^2 + a |u|^2 dx + \intlog -\intR F(u) dx.
\end{equation}
But note that $I$ is not well define in $\Hum$. So, following the idea introduced by Stubbe \cite{[21]}, we will consider the slightly smaller Hilbert space
\beq \label{X}
X = \left\{ u\in \Hum ; \intR \ln(1+|x|)u^2(x) dx < \infty \right\} \subset \Hum ,
\eeq
where the part of $I$ involving the $ \ln |\cdot| $ function will be finite. Also, we need to guarantee that $ \intR F(u) dx < + \infty $. From Moser-Trudinger inequality, Lemma \ref{l23}, for $ \alpha= 4\pi \tau > 0 $, with $\tau >1,$ we have
\begin{align*}
\intR F(u) dx & \leq \intR [|u|^2 + b_2 |u|^p(e^{\alpha u^2}-1)] dx \\
& \leq  ||u||_{2}^{2} + b_2 ||u||_{r_1 p}^{p}\left(\intR (e^{r_2 \alpha u^2}-1) dx\right)^{\frac{1}{r_2}} < + \infty ,  
\end{align*}
for all $ u\in X $ with $ \frac{1}{r_1}+\frac{1}{r_2}=1 $, with $r_2 \sim 1. $

Therefore, $ I: X \RA \R $ as given in (\ref{I}) is well-defined. As we will see, $ I $ is of class $ C^1 $ with Gateux derivative given by 
\beq \label{I'}
I'(u)(\varphi)=\intR [\nabla u \nabla \varphi + a u \varphi ] dx + \intLog - \intR f(u) \varphi dx  \ ,
\eeq
for all $ u \in X $, $ \varphi\in \Hum $.

Concerning to the Hilbert space $ X $, define
\beq\label{ast}
||u||^{2}_{\ast}=\intR \ln(1+|x|) u^2(x) dx, \ \forall \ u \in X.
\eeq
If $ ||u||^2 = ||\nabla u||_{2}^{2}+||u||_{2}^{2} $ is the usual norm of $ \Hum $, the expression $$ ||u||_{X}^{2} = ||u||^{2}+||u||_{\ast}^{2}$$ is a norm on $X$.

Following the ideas of some recent papers, as for example \cite{[guo],[7], [wen] , [8]}, we can define the Pohozaev's functional associated with  (\ref{P}), $ P:X \RA \R $, by
\begin{align*}P(u)& = a\intR u^2(x)dx +\intlog + \dfrac{1}{4}( \intR u^2(x) dx )^2 - 2\intR F(u) dx
\end{align*}
and based on that, a key auxiliar functional  $ J:X\RA \R $, given by 
\begin{align*}
J(u)&=\dfrac{5}{2}I'(u)(u)-P(u)\\
&=\dfrac{5}{2}\intR |\nabla u|^2 dx + \dfrac{3}{2}\intR a u^2 dx + \dfrac{3}{2}\intlog \\ &- \dfrac{1}{4}\left( \intR u^2(x) dx \right)^2 
 + 2\intR F(u) dx - \dfrac{5}{2}\intR f(u)u dx .
\end{align*}
Then, similarly as in  \cite{[8] , [6] , [10]}, we can verify that any weak solution $ u $ of (\ref{P}) satisfies $ P(u)=0 $. To be precise, we have the following lemma

\bl\label{l26}
Suppose that $u\in X$ is a weak solution to (\ref{P}). Then, 
$P(u)=0,$
i.e., the Pohozaev functional is zero over all  the  weak solutions of (\ref{P}).
\el

Our strategy is to find a sequence $(\un)\in X$ satisfying 
\beq\label{i4}
I(\un)\RA c_{mp} \ , \ ||I'(\un)||_{X'}(1+||\un||_X)\RA 0 \  \mbox{and} \ J(\un)\RA 0 ,
\eeq
where
\beq\label{cmp}
c_{mp}=\inf\limits_{\gamma \in \Gamma}\max\limits_{t\in [0, 1]}I(\gamma(t)),
\eeq
with $ \Gamma = \{ \gamma \in C([0, 1], X) \ ; \ \gamma(0)=0 , I(\gamma(1))<0 \} $, is the mountain pass level for $ I $. 

In order to verify that such sequence is bounded in $\Hum$, we will need the following condition 
$$ \mbox{ there exists }\ \theta > \dfrac{16}{5} \ \mbox{ such that}\   f(s)s \geq \theta F(s) > 0, \ \mbox{for all }\  s> 0. \leqno{(f_3)}$$
Here it is interesting point out that some usual conditions applied to nonlinearities of the form $ f(u)=b|u|^{p-2}u $, for $ b>0 $, $ p>\theta $, does not help in our case, once Lemma \ref{l23} demands that we already know that $ ||\nabla \cdot ||_{2}^{2}\leq 1 $ to obtain the boundedness for the exponential integral. \\
As we pointed above, some kind of stronger condition is needed to verify that the integral involving the exponential is bounded, in particular, to our wished sequence (\ref{i5}). So, we include the following  
$$ \mbox{ there  exist}\  q>4 \ \mbox{ and}\  C_q> \dfrac{[6(q-2)]^{\frac{q-2}{2}}}{q^{\frac{q}{2}}}\dfrac{S_{q}^{q}}{\rho^{q-2}} \ \mbox{ such that}\  F(s) \geq C_q |s|^q , \ \mbox{for all}\  s\geq 0 , \leqno{(f_4)}$$ 
for $ S_q, \rho >0 $ to be defined in Lemma \ref{l211}. Note that, in particular, inequalities (\ref{eq1}) and (\ref{eq2}) hold for $ p=q $, given as above. So, along of all paper, we will consider $ p=q $. A prototype of a nonlinearity satisfying condition $ (f_1)-(f_4) $ is given by $ f(s)=0 $, for $ s\leq 0 $, and  $$ f(s)=C_q \left\{ \begin{array}{ll}
s^{q-1} \ \ \ \ \ \ \ \ \ \ \ \ \ \ \ \ , \ \mbox{if} \ 0 \leq s \leq 1 \\
s^{q-1} e^{4\pi (s^2 - 1)} \ , \ \mbox{if} \ 1< s , 
\end{array} \right. $$ for $ C_q>0 $ sufficiently large and $ q>4 $.

Now we are ready to enunciate our first main result, concerned with the existence of a solution for (\ref{P}) at the mountain pass level and being a ground state solution.

\bt \label{t11}
Suppose that $ f $ satisfies $ (f_1)-(f_4) $, $ \alpha_0 = 4\pi $, $ q>4 $ and $ C_q>0 $ is sufficiently large. Then, we have the following: 
\begin{itemize}
\item[(i)] $ c_{mp}>0 $ and (\ref{P}) has a solution $ u\in X\setminus \{0\} $ with $ I(u)=c_{mp} $. 
\item [(ii)] (\ref{P}) has a ground state solution, i.e., a solution $ u\in X\setminus \{0\} $ such that $ I(u)=c_g = \inf\{ I(v) \ ; \ v\in X\setminus \{0\} \textrm{ \ is a solution of (\ref{P})}\} $. 
\end{itemize}
\et

Our second main result follows the ideas of \cite{[10] , [6]} and deals with a symmetric setting with respect to a suitably action of subgroups of $ O(2) $. 

First we need to introduce some notations. Let $ G $ be a subgroup of the orthogonal group $ O(2) $ and $ \tau : G \RA \{ -1 , 1\} $ a group homomorphism. Then, $ (G, \tau) $ defines a group action of $ G $ on $ X $ by 
$$
[A \ast u](x) \doteq \tau(A)u(A^{-1}x) \textrm{ \ , for \ } A\in G, u\in X, x\in \Rdois.
$$
In this point of view, we have the following invariant space
$$
X_G \doteq \{u\in X \ ; \ A\ast u = u \textrm{ \ for all \ } A\in G\}\subset X .
$$

Also, in view of getting the invariance of $I$ under the group action, we will need to change the condition $ (f_1) $ for 
 $$\left\{ \begin{array}{c}f: \R \RA \R \ \mbox{is continuous, odd and has critical exponential growth with } \ \alpha_0 = 4 \pi,   \\
 |f(s)|\leq C e^{4\pi s^2}, \forall \ s \in \R
 \end{array}\right. \leqno{(f_1')}$$
 
Thus, we are ready to establish a multiplicity result.

\bt\label{t12}
Suppose that $ f $ satisfies $ (f_1')-(f_4) $, $ \alpha_0 = 4\pi $, $ q>4 $ and $ C_q>0 $ is sufficiently large. Let $ G, \tau $ be as above, and assume that $ X_G \neq \{0\} $. Define 
\beq\label{eq30}
c_{mp, G}=\inf\limits_{\gamma \in \Gamma_{G}}\max\limits_{t\in [0, 1]}I(\gamma(t)),
\eeq
with $ \Gamma_G = \{ \gamma \in C([0, 1], X_G) \ ; \ \gamma(0)=0 , I(\gamma(1))<0 \} $. Then, we can prove the following 
\begin{itemize}
\item [(i)] $ c_{mp, G}>0 $ and equation (\ref{P}) has a solution $ u\in X_G \setminus \{0\} $ with $ I(u)= c_{mp, G} $. 
\item [(ii)] Equation (\ref{P}) has a ground state solution $ u\in  X_G \setminus \{0\} $ that is $ G- $invariant, i.e. $ u\in X_G $ satisfies (\ref{P}) and $ I(u)=c_{g, G}=\inf\{I(v) \ ; \ v\in X_G \setminus \{0\}  \textrm{ \ is a solution of (\ref{P})} \}  $.
\end{itemize}
\et

\bo\label{obs10}
The above theorem tell us that problem (\ref{P}) has infinitely many nonradial sign-changing solutions. If the homomorphism $\tau$ is nontrivial and $ A\in G $ satisfies $ \tau(A)=-1 $, then, for every $ u\in X_G $, we have $ u(x)=[A\ast u](x)=\tau(A)u(A^{-1}x)=-u(A^{-1}x) $. Hence, if $ x\in W=\{ y\in \Rdois \ ; \ Ay=y\}  $, then $ A^{-1}x=x $ and $ u(x)=-u(x) $, which implies that $ u$ vanishes on the set $ W $. On the other side, if $ u\neq 0 $, $ u $ changes sign.  
\eo

We also would like to point out, as the authors in \cite{[10]}, that no analogue of Theorem \ref{t12} has been proved yet for higher dimensions, that is, in the case $ N\geq 3 $ and $ \lambda >0 $. Next, we present an example as an illustration for this setting. For more examples and further discussion we refeer to \cite{[10]}. 

\be
Let $ G=O(2) $ and $ \tau \equiv 1 $. Then, $ u(x)=[A\ast u](x)=u(A^{-1}x) $, for all $ u\in X_G $ and $ A\in G $. Remember that a function $ g $ is radial if and only if $ g(\rho x) = g(x) $, for all rotations $ \rho $ and $ x\in \Rdois $. Therefore, in this case $ X_G $ consists of all radial functions in $ X $.
\ee

Throughout the paper, we will use the following notations: $ \Ls $ denotes the usual Lebesgue space with norm $ ||\cdot ||_s $ \ ; \ $ X' $ denotes the dual space of $ X $ \ ; \ $ B_r(x) $ is the ball centred in $ x $ with radius $ r>0 $ \ ;  \ $ C, C_1, C_2, ... $ will denote different positive constants whose exact values are not essential to the exposition of arguments. 

The paper is organized as follows: in section 2 we present some technical and essential results, some of them already derived in previous works and whose application to our problem is immediate. Section 3 consists in the proof of a key proposition and our first main result. Finally, in section 4, we give some details about the symmetric setting and prove the second main result.

\section{Preliminary Results}

In this section, we present some additional framework, observations and provide technical lemmas that will help in the proofs of our main results. We start pointing that, from the fact that $ a>0 $, we can endow $ \Hum $ with the  scalar product  given by
\begin{equation*}
 \langle u , v \rangle = 
\intR [\nabla u \nabla v + a u v ] dx , \ \textrm{ for \ } u, v \in \Hum ,
\end{equation*}
with the corresponding norm  $\|u\|=  \sqrt{\langle u , u\rangle },$
which one is equivalent to the usual norm.  
 So,  we will work with both without distintion.

Next, inspired by \cite{[6]}, we define three auxiliar symmetric bilinear forms 
\beqa
(u , v) \mapsto B_1(u, v)=\intR \intR \ln(1+|x-y|)u(x)v(y) dx dy ,
\eeqa
\beqa
(u , v) \mapsto B_2(u, v)=\intR \intR \ln\left(1+\dfrac{1}{|x-y|}\right)u(x)v(y) dx dy ,
\eeqa
\beqa
(u , v) \mapsto B_0(u, v)=B_1(u, v)-B_2(u, v)=\intR \intR \ln(|x-y|)u(x)v(y) dx dy. 
\eeqa
These definitions are understood to being over measurable function $u, v: \Rdois \RA \R $, such that the integrals are defined in the Lebesgue sense. Then, since $ 0 \leq \ln(1+r)\leq r $, for $ r>0 $, and by Hardy-Littlewood-Sobolev Inequality (HLS) \cite{[15]}, for $ u, v \in L^{\frac{4}{3}}(\Rdois) $,
\begin{align} \label{b2}
|B_2(u, v)| 
& \leq \intR \intR \dfrac{1}{|x-y|}u(x)v(y) dx dy 
 \leq K_0 ||u||_{\frac{4}{3}}||v||_{\frac{4}{3}} , 
\end{align}
where $ K_{0}>0 $ is the HLS constant. 

Next, we define the functionals  $V_1: \Hum \RA [0, \infty],$ $V_2: L^{\frac{8}{3}}(\Rdois) \RA [0, \infty) $ and $V_0: \Hum \RA \R \cup \{\infty\},$ given by 
$V_1(u)=B_1(u^2 , u^2), $ $V_2(u)=B_2(u^2 , u^2) $ and $V_0(u)=B_0(u^2 , u^2),$ respectively.

As a consequence of (\ref{b2}), we have
\beq \label{v2}
|V_2(u)|\leq K_0||u||_{\frac{8}{3}}^{4} \ , \ \ \forall \ u\in L^{\frac{8}{3}}(\Rdois) ,
\eeq
so $ V_2 $ takes finite values over $ L^{\frac{8}{3}}(\Rdois)\subset \Hum $. Also, observing that
\beq \label{2.3}
\ln(1+|x-y|) \leq \ln(1+|x|+|y|) \leq \ln (1+|x|)+ \ln(1+|y|), \textrm{ \ for \ } x, v\in \Rdois,
\eeq

we can estimate, applying Hölder Inequality,
\beq\label{b1}
B_1(uv, wz) \leq ||u||_{\ast}||v||_{\ast}||w||_2 ||z||_2 + ||u||_2 ||v||_2 ||w||_{\ast}||z||_{\ast} ,
\eeq
for all $ u, v, w, z \in \Ldois $. \\

We need the following technical lemmas proved in \cite{[6]}, so we will omit their proofs here.

\bl \label{l21}
(\cite{[6]}, Lemma 2.1) Let $(\un)$ be a sequence in $\Ldois $ and $ u\in \Ldois \setminus \{0\} $ such that $\un \RA u $ pointwise a.e. on $\Rdois $. Moreover, let $(\vn)$ be a bounded sequence in $\Ldois$ such that
$$
\sup\limits_{n\in \N}B_1(\un^2, \vn^2) < \infty .
$$
Then, there exist $n_0\in \N$ and $ C>0 $ such that $ ||\un||_{\ast}< C $, for $ n\geq n_0 $. If, moreover, 
$$
B_1(\un^2, \vn^2)\RA 0 \textrm{ \ \ and \ \ } ||\vn||_{2}\RA 0, \textrm{ \ as \ } n\RA \infty ,
$$
then 
$$
||\vn||_{\ast}\RA 0 \textrm{ \ ,  as \ } n \RA \infty .
$$
\el

\bl\label{l212}
(\cite{[6]}, Lemma 2.6) Let $ (\un) $, $ (\vn) $ and $ (\wn) $ be bounded sequences in $ X $ such that $ \un \rightharpoonup u $ in $ X $. Then, for every $ z\in X $, we have $ B_1 (\vn \wn \ , \ z(\un - u))\RA 0 $, as $ n\RA + \infty $.
\el

\bl \label{l22} (\cite{[6]}, Lemma 2.2)
(i) The space $ X $ is compactly embedded in $ \Ls $, for all $ s\in [2, \infty) $. \\
(ii) The functionals $ V_0, V_1, V_2 $ and $ I $ are of class $ C^1 $ on $ X $. Moreover, $ V_{i}'(u)(v)=4B_i(u^2, uv) $, for $ u, v \in X $ and $ i=0, 1, 2 $. \\
(iii) $ V_2 $ is continuous (in fact continuously differentiable) on $\Loito$ .\\
\el

Moreover, in view of Moser-Trudinger inequality (see also \cite{ [9],Lam}), we can prove the following very usefull result.
 




\bl\label{l23}
\cite{[5]} If $\alpha >0$ and $ u\in \Hum $, then
$$
\intR \left(e^{\alpha |u|^2} - 1 dx \right) < \infty .
$$
Moreover, if $ ||u||^2\leq 1 $, $ ||u||_{2}\leq M < \infty $ and $ \alpha < \alpha_0 = 4 \pi $, then there exists $ C_{\alpha}=C(M, \alpha) $, such that
$$
\intR \left(e^{\alpha |u|^2} - 1 \right) dx < C_\alpha .
$$
\el

Therefore, from above lemmas, letting $ p>2 $, $\eta\in (0,1) $, $ \alpha_0 = 4\pi ,$ $\tau >1,$ fixed, and $ ||u||\leq \delta $, under Lemma \ref{l23} conditions, one can easily verify that there exists constants  $C_1 , C_2 >0$ such that
\beq \label{eq5}
\intR F(u) dx \leq \dfrac{\eta}{2} ||u||^{2} + C_1 ||u||^{p}.
\eeq
and
\beq \label{eq6}
\intR f(u)u \ dx \leq \intR |f(u)u| dx \leq \eta ||u||^{2} + C_2 ||u||^{p} ,
\eeq

The next two lemmas will show that the functional $ I $ has a mountain pass geometry.

\bl \label{l25}
There exists $ \rho > 0  $ such that
\beq \label{eq7}
m_\beta = \inf\{I(u) ; u\in X \ , \ ||u||=\beta \}>0 \ , \textrm{ \ for \ } 0 < \beta \leq \rho
\eeq
and 
\beq \label{eq8}
n_\beta = \inf\{I'(u)(u) ; u\in X \ , \ ||u||=\beta \}>0 \ , \textrm{ \ for \ } 0 < \beta \leq \rho .
\eeq
\el

\textbf{Proof:} Let $ u\in X \setminus \{0\} $ such that $ ||u||\leq \delta  $ and $ \eta, \tau $ as in equations (\ref{eq5}) and (\ref{eq6}). Then, by Sobolev embeddings and equation (\ref{eq5}),
$$
I(u) = \dfrac{1}{2}||u||^{2} + \dfrac{V_1(u)}{4} - \dfrac{V_2(u)}{4} - \intR F(u) dx \geq \left(\dfrac{1}{2}-\frac{\eta}{2}\right)||u||^{2} [ 1- C_3||u||^2 -C_4 ||u||^{p-2}] .
$$
Consequently, for $ \rho > 0 $ sufficiently small, (\ref{eq7}) is valid. Similarly, from (\ref{eq6})
\begin{align*}
I'(u)(u) 
& \geq \left(1-\eta \right)||u||^2[1-C_2||u||^2 -C_3 ||u||^{p-2}].
\end{align*}
Therefore, for $ \rho > 0 $ sufficiently small, we have (\ref{eq8}). 
\cqd
\bo \label{obs4}
Lemma \ref{l25} tell us that there exists $ \rho >0 $ and $ b_0 >0 $ such that $I(u) > 0 $, for all $u\in X$ with $ 0< ||u||\leq \rho $ and $ I(u)>b_0 >0 $, for all $ u\in X $ such that $ ||u||=\rho $.
\eo

\bl\label{l27}
Let $ u\in X\setminus \{0\} $, $ t>0 $ and $ q> 4 $. Then,
$$
\lim\limits_{t\RA 0} I(tu) = 0 \ \ , \ \ \sup\limits_{t>0} I(tu) < +\infty \ \ \mbox{and} \ \ I(tu)\RA + \infty \ , \ \mbox{as} \ t \RA + \infty .
$$
\el
\textbf{Proof:} Let $ u\in X\setminus \{0\} $. First of all, observe that
$$
I(tu)= \dfrac{t^2}{2}||\nabla||_{2}^{2}+\dfrac{at^2}{2}||u||_{2}^{2}+\dfrac{t^{4}}{4}V_0(u)-\intR F(tu) dx .
$$
From $ (f_4) $, since $ q>4 $,
$$
I(tu) \leq \dfrac{t^2}{2}||\nabla||_{2}^{2}+\dfrac{at^2}{2}||u||_{2}^{2}+\dfrac{t^{4}}{4}V_0(u)-C_q t^{q}||u||_{q}^{q} \RA - \infty ,
$$
as $ t\RA + \infty $. Now, observe that, for $ \frac{1}{r}+\frac{1}{r'}=1  $, with $ r\sim 1 $, a fixed $ u\in X\setminus \{0\} $ and $ t>0 $ sufficiently small, $ r ||t u||^2 < 4\pi $, and then
\beq \label{eq13}
\left| \intR F(tu) \right| \leq \dfrac{t^2}{2}||u||_{2}^{2}+b_2 t^{q}||u||_{r' q}^{q}\left(\intR \left(e^{r (t^{\frac{5}{2}}u)^2} - 1 \right) dx \right)^{\frac{1}{r}} \leq \dfrac{t^2}{2}||u||_{2}^{2} +C_1 t^{q}||u||_{r' q}^{q} .
\eeq
From (\ref{eq10}) and (\ref{eq13}), $ I(tu)\RA 0 $ as $ t\RA 0 $. Therefore, from the two assertions proved, also follows that $\sup\limits_{t>0}I(tu)< \infty  $. Finally, from these two assertions and the fact that $ I $ is $ C^1 $, follow that $ \sup\limits_{t>0} I(tu) < +\infty $. 

\cqd

One can easily see, applying the Intermediate Value Theorem, that $ 0<b_0 \leq c_{mp}<+\infty $. Next, we present a necessary proposition to obtain our desired sequence satisfying (\ref{i4}).

\bp \label{l210}
(\cite{[13]}, Proposition 2.8) Let $ X $ be a Banach Space, $ M_0 $ be a closed subspace of the metric space $ M $ and $ \Gamma_0 \subset C(M_0 , X)$. Define $ \Gamma = \{ \gamma \in C(M, X) \ ; \ \gamma\vert_{M_0} \in \Gamma_0\} $. If $ \varphi \in C^{1}(X, \R) $ satisfies 
$$
\infty > c \doteq \inf\limits_{\gamma \in \Gamma}\sup\limits_{u\in M} \varphi(\gamma(u)) > a \doteq \sup\limits_{\gamma_0 \in \Gamma_0}\sup\limits_{u\in M_0} \varphi(\gamma_0 (u)) ,
$$
then, for every $ \varepsilon \in \left(0 , \frac{c-a}{2}\right) $ and $ \gamma \in \Gamma $ with $ \sup\limits_{u\in M}\varphi(\gamma(u)) \leq c + \varepsilon $, there exists $ u\in X $ such that 
\begin{itemize}
\item [ \textbf{(a)} ]$ c-2\varepsilon \leq \varphi(u) \leq c + 2 \varepsilon $ , 

\item [\textbf{(b)} ] $ (1+||u||_X)||\varphi'(u)||_{X'}\leq \dfrac{8\varepsilon}{\delta} $ , 

\item [\textbf{(c)} ]$ dist(u , \gamma(M)) \leq 2 \delta $.
\end{itemize}
\ep

As mentioned in the works of Cingolani -Jeanjean \cite{[cjj]} and Du -Weth \cite{[10]}, the study of the planar case of Choquard equation has a major difficulty to  ensure  when a Cerami sequence will present a strong convergent subsequence in $ X $. All in all, the next lemma is a key tool in our study, once it provides a compactness condition.  

\bl\label{l29}
There exists a sequence $(\un)\subset X$ satisfying
\beq\label{eq15}
I(\un)\RA c_{mp} \ , \ ||I'(\un)||_{X'}(1+||\un||_{X})\RA 0 \textrm{ \ \ and \ \ } J(\un) \RA 0 \ , \ \textrm{ \ as \ } n\RA + \infty .
\eeq
\el
\textbf{Proof:} Consider the Banach space $ \tilde{X}=\R \times X $ equipped with the norm $ ||(s, v)||_{\tilde{X}}=(|s|^{2}+||v||_X^{2} )^{\frac{1}{2}} $. Let's define $ \varphi : \tilde{X}\RA X $ by $ \varphi(s, v)(x) = e^{\frac{5}{2}s}v(e^s x) $, for $ s\in \R $, $ v\in X $ and $ x\in \Rdois $. Note that $ \varphi $ is continuous. Define also $ \phi : \tilde{X}\RA \R $ setting $ \phi (s, v) =I(\varphi(s, v))  $, for $ s\in \R $ and $ v\in X $. So, by a change of variables, one can easily see that  
\begin{eqnarray} \label{eq16}
\lefteqn{\phi (s, v) = }\\ && \dfrac{e^{5s}}{2} \intR |\nabla v|^{2} dx + \dfrac{a}{2}e^{3s} \intR |v|^2 dx + \dfrac{e^{6s}}{4}V_0(v) - \dfrac{se^{6s}}{4}\left(\intR |v|^2 dx \right)^2 - \dfrac{1}{e^{2s}}\intR F(e^{\frac{5}{2}s}v) dx .\nonumber
\end{eqnarray}
Clearly, $ \phi \in C^{1}(\tilde{X}, \R) $ and we have
\begin{align}
\partial_s \phi(s, v) & = \dfrac{5}{2}e^{5s} \intR |\nabla v|^{2} dx + \dfrac{3a}{2}e^{3s} \intR |v|^2 dx + \dfrac{3}{2}e^{6s} V_0(v) - \dfrac{3}{2}se^{6s}\left(\intR |v|^2 dx \right)^2 \nonumber \\
& - \dfrac{e^{6s}}{4}\left(\intR |v|^2 dx \right)^2 +\dfrac{2}{e^{2s}}\intR F(e^{\frac{5}{2}s}v) dx - \dfrac{5}{2}e^{\frac{s}{2}}\intR f(e^{\frac{5}{2}s}v) v dx . \label{eq17}
\end{align}
Consequently, from (\ref{eq17}), 
\beq \label{eq19}
\partial_s \phi(s, v) = J(\varphi(s, v)) \textrm{ \ , for \ }  (s , v) \in \tilde{X} .
\eeq 
Moreover, since $ v \mapsto \varphi(s, v) $ is linear for every fixed $ s\in \R $, we have 
\beq \label{eq20}
\partial_v \phi(s, v)(w) = I'( \varphi(s, v))\varphi(s, w) \textrm{ \ , for \ } s\in \R, v,w \in X .
\eeq
Define the minimax value $ \tilde{c}_{mp} $ for $ \phi $ by
$$
\tilde{c}_{mp}=\inf\limits_{\tilde{\gamma} \in \tilde{\Gamma}}\max\limits_{t\in [0, 1]}\phi(\tilde{\gamma}(t)) ,
$$
where $ \tilde{\Gamma}=\{ \tilde{\gamma}\in C([0, 1], \tilde{X}) \ ; \ \tilde{\gamma}(0)=0 , \phi(\tilde{\gamma}(1))<0\} $. Defining $\tilde{\gamma}$ by $ \tilde{\gamma}(s)=(0, \gamma(s)) $, one can easily verify that $ \Gamma = \{ \varphi \circ \tilde{\gamma} \ ; \ \tilde{\gamma} \in \tilde{\Gamma}\} $. So, we conclude that the minimax values of $ I $ and $ \phi $ coincide, i.e., $ c_{mp}=\tilde{c}_{mp} $. Now, from the definition of $ c_{mp} $, for $ n\in \N $, there exists $ \gamma_n \in \Gamma $ such that
\beq\label{eq21}
\max\limits_{t\in [0, 1]}\phi(0, \gamma_n (t)) = \max\limits_{t\in [0, 1]} I(\gamma_n(t))\leq c_{mp}+\dfrac{1}{n^2} .
\eeq
Let $ X=\tilde{X} $, $ M_0 = \{ 0, 1\} $, $M=[0, 1]$, $ \Gamma = \tilde{\Gamma} $ and $ \varphi = \phi $ in Lemma \ref{l210} and, let also $ \varepsilon_n = \frac{1}{n^2} $, $ \delta_n = \frac{1}{n} $ and $ \tilde{\gamma}_n =(0, \gamma_n (t)) $. From $ c_{mp}\geq b_0 >0 $, $\varepsilon_n \in (0, \frac{c}{2})$, for $ n $ sufficiently large, and Lemma \ref{l210}, we obtain a sequence $ (\sn , \vn)\in \tilde{X} $ such that, as $ n\RA + \infty $, we have 

\textbf{(a)} $ \phi(\sn, \vn ) \RA c_{mp} $ , 

\textbf{(b)} $ ||\phi'(\sn , \vn)||_{\tilde{X}'}(1+||(\sn , \vn)||_{\tilde{X}})\RA 0 $ , 

\textbf{(c)} $ dist((\sn , \vn ), \{0\}\times \gamma_n ([0, 1]))\RA 0 $ . 

\noindent From (c), we conclude that $ \sn \RA 0 $. Observe that, from (\ref{eq19}) and (\ref{eq20}), for $ (h, w)\in \tilde{X} $,
\beq\label{eq22}
\phi'(\sn , \vn)(h, w)=I'(\varphi(\sn , \vn))\varphi(\sn , w)+J(\varphi(\sn , \vn))h .
\eeq
So, taking $ h=1 $ and $ w=0 $ in (\ref{eq22}),
\beq\label{eq23}
\phi'(\sn , \vn)(1, 0)=I'(\varphi(\sn , \vn))\varphi(\sn , 0)+J(\varphi(\sn , \vn))1=J(\varphi(\sn , \vn)).
\eeq
From (b) and (\ref{eq23}), $ J(\varphi(\sn , \vn))\RA 0 $. Thus, setting $ \un = \varphi(\sn , \vn $), from (a), follows that
$$
I'(\un)\RA c_{mp} \textrm{ \ \ and \ \ } J(\un)\RA 0 \ \ , \ \textrm{ as \ } n \RA + \infty .
$$
To conclude, given $v\in X $, consider $ \wn (x) = e^{-\frac{5}{2}\sn}v(e^{-\sn}x) $, for $ x\in \Rdois $. Then, taking $ h=0 $ and $ w=\wn $ in (\ref{eq22}), we have
$$
(1+||\un||_{X})|I'(\un)v|\leq (1+ ||\un||_{X})||\phi'(\sn , \vn) ||_{\tilde{X}'}||\wn||_{X}= o(1)||\wn||_{X}.
$$
But, as one can easily verify, as $ n\RA + \infty $ with $ o(1)\RA 0 $ uniformly for $ v\in X $,
$$
||\wn||_{X}^{2}= e^{-3\sn} \intR |\nabla v|^2 dx + e^{-5\sn} \intR [a|v|^2 + \ln(1+e^{\sn}|x|) v^2] dx = (1+o(1))||v||_{X}^{2} .
$$
From these last two estimates, we have
$$
(1+||\un||_{X})|I'(\un)v|= o(1)||\wn||_{X}=o(1)(1+o(1))||v||_{X}^{2}\RA 0 \ \ , \ as \ n\RA + \infty .
$$
Therefore, $ ||I'(\un)||_{X'}(1+||\un||_{X})\RA 0 $. 

\cqd

\bp\label{p21}
Let $(\un) \subset X$ satifying (\ref{eq15}). Then, $ (\un) $ is bounded in $ \Hum $.
\ep
\noindent  \textbf{Proof:} From (\ref{eq15}) and $ (f_3) $,
\beq\label{eq24}
c_{mp} + o(1)  = I(\un)-\dfrac{1}{6}J(\un) \geq \dfrac{1}{12}||\un||^{2}+ \intR \left(\dfrac{5}{12}\theta - \dfrac{4}{3} \right) F(\un) dx \geq \dfrac{1}{12}||\un||^{2} .
\eeq
Therefore, $ (\un) $ is bounded in $ \Hum $. 

\cqd

\bl\label{l211}
Let $ (\un) \subset X $ satisfying (\ref{eq15}), $ q>4 $ and $ C_q >0 $ sufficiently large. Then, there exists $ \rho_0 >0$ sufficiently small such that
$$
\limsup\limits_{n} ||\un|| < \rho_0 .
$$
\el
\textbf{Proof:} From (\ref{eq24}), $ 12c_{mp} + o(1)\geq ||\un||^{2} $ and so $ \limsup\limits_{n} ||\un||^2  \leq 12 c_{mp}$.

Consider the set $ A=\{u\in X\setminus \{0\} \ ; \ V_0(u)\leq 0 \} $. Lets prove that $ A\neq \emptyset $. In fact, let $ u\in X\setminus \{0\} $ and define $ u_t(x)=t^{\frac{5}{2}}u(tx) $. Then, $ V_0(u_t)=t^6 V_0(u)-t^6 \ln t ||u||_{2}^{4} $. Thus, for $ T>0 $ sufficiently large, $ V_0(u_T)<0 $. Therefore, $ u_T \in A $ and $ A\neq \emptyset $. 

Also, by Sobolev embeddings, there exists a constant $ C>0 $ such that $ ||u||\geq C||u||_{q} $. So, it makes sense to define
$$
S_q(v)=\dfrac{||v||}{||v||_{q}} \textrm{ \ \ \ and \ \ \ } S_q = \inf\limits_{v\in A} S_q(v) \geq \inf\limits_{v\in \Hum\setminus \{0\} } S_q(v) >0 .
$$
Hence, we are ready to find an estimative for $ c_{mp} $. For $ \psi \in A $, we see that 
\begin{eqnarray} \label{eq25}
c_{mp} \leq \max\limits_{t\geq 0} \left\{\dfrac{t^2}{2}||\psi||^2- C_q t^q ||\psi||_{q}^{q}\right\} &\leq& \max\limits_{t\geq 0} \left\{\dfrac{S_q(\psi)^2}{2}t^2||\psi||_{q}^2- C_q t^q ||\psi||_{q}^{q}\right\} \nonumber \\
& \leq & \left(\dfrac{1}{2}-\dfrac{1}{q}\right) \dfrac{S_q(\psi)^{\frac{2q}{q-2}}}{(qC_q)^{\frac{2}{q-2}}} \nonumber
\end{eqnarray}
Taking the infimum over $ \psi \in A $, we obtain
$$
\limsup\limits_{n} ||\un||^2 \leq 6\dfrac{(q-2)}{q} \dfrac{S_q^{\frac{2q}{q-2}}}{(qC_q)^{\frac{2}{q-2}}} \leq \rho^2 .
$$

\cqd

\section{Proof of Theorem 1.1}

In this section we will finish the proof of Theorem 1.1. In the following, we define a function $ y \ast u : \Rdois \RA \R $ by $ [y \ast u ] (x) = u(x-y) $, for $ u\in X $ and $ y, x \in \Rdois $. The next proposition will provide us with sufficient condition to find a nontrivial critical point for $ I $ in $ X $.

\bp\label{p31}
Let $ q>4 $ and $ (\un)\subset X $ satisfying (\ref{eq15}). Then, after passing to a subsequence, one between  two alternatives occurs:

\textbf{(a)} $ ||\un||\RA 0 $ and $ I(\un)\RA 0 $, as $ n\RA + \infty $. 

\textbf{(b)} There exists points $ y_n \in \mathbb{Z}^2 $ such that $ y_n \ast \un \RA u $ in $ X $, for some critical point $ u\in X\setminus \{0\} $ of $ I $.
\ep
\textbf{Proof:} From Proposition \ref{p21} and Lemma \ref{l211}, $ (\un) $ is bounded in $ \Hum $ and $ ||\un|| $ is sufficiently small, for all $ n\in \N $. Suppose that $ (a) $ does not occur for any subsequence of $ (\un) $. \\
\textbf{Claim 1:} $ \liminf\limits_{n\RA + \infty} \sup\limits_{y\in \mathbb{Z}^2} \ds_{B_2 (y)} \un^{2}(x)dx >0 $. 

Suppose the contrary happens. So, by Lion's Lemma \cite{[dg]}, after passing to a subsequence, satisfies $ \un \RA 0 $ in $\Ls $, for $ s > 2 $. From this and (\ref{eq2}), we have
$$
||\un||^{2}+V_1 (\un)=I'(\un)(\un) +V_2 (\un) +\intR f(\un)\un \RA 0 \ \ , \textrm{ \ as \ } n\RA + \infty .
$$
Consequently, since $ ||\un||\geq 0 $ and $ V_1 (\un) \geq 0 $, $ ||\un||\RA 0 $ and $  V_1 (\un)\RA 0 $. So, we conclude that
$$
I(\un)=\dfrac{1}{2}||\un||^2 + \dfrac{1}{4}[V_1(\un)-V_2(\un)] - \intR F(\un) dx \RA 0 .
$$
But this contradicts the supposition that (a) does not occur. Therefore, the claim is valid. 

Hence, passing to a subsequence, if necessary, there exists $ (y_n)\subset \mathbb{Z}^2 $ such that $ \until =y_n \ast \un \in X $, for all $ n \in \N $. Once $ (\un) $ is bounded in $ \Hum $, so is $ (\until) $ and then $ \until \rightharpoonup u $ in $ \Hum\setminus\{0\} $. Thus, we can assume that $ \until(x)\RA u(x) $ a.e. in $ \Rdois $. \\

Therefore, from the boundedness of $(\un)$  in $ \Hum $, (\ref{eq15}) and (\ref{eq6}), we deduce that
$$
B_1 (\until^2 , \until^2 ) = V_1 (\until) = V_1 (\un) = o(1) + V_2 (\un) + \intR f(\un) \un dx -||\un||^2 , 
$$
remains bounded  for all $ n $, i.e., $ \sup\limits_{n} B_1 (\until^2 , \until^2 ) < + \infty $. So, applying Lemma \ref{l21}, $ ||\until||_{\ast} \leq C $, for a constant $ C>0 $. As $ (\until) $ is already bounded in $ \Hum $, it follows that $ (\until) $ is bounded in $ X $. From reflexiveness of $ X $, passing to a subsequence if necessary, $ \until \rightharpoonup u $ in $ X $. For this reason, $ u\in X $ and from Lemma  \ref{l22}- (ii),  we  get $ \until \RA u $ in $ \Ls $, for $ s\geq 2 $. \\
\textbf{Claim 2:} $ I'(\until)(\until - u)\RA 0 $, as $ n\RA + \infty $. \\
In fact, for all $n \in \N,$
\beq\label{eq28}
|I'(\until)(\until - u)|=|I'(\un)(\un - (-\yn) \ast u)| \leq ||I'(\un)||_{X'}(||\un||_X + ||(-\yn) \ast u||_X).
\eeq
Then, we first seek for an useful inequality for $ ||(-\yn) \ast u||_X $. If $|\yn|\RA +\infty$, then, for $ x\in \Rdois $,
$$
\ln(1+|x-\yn|)-\ln(1+|\yn|)=\ln \left(\dfrac{1+|x-\yn|}{1+|\yn|}\right) \RA 0 , n\RA + \infty .
$$
Therefore, there exists $ C_1 > 0 $ such that $ \ln(1+|x-\yn|)\geq C_1 \ln(1+|\yn|) $. \\
Now, suppose that $ (\yn)\subset \mathbb{Z}^2 $ converges to $ y_0 \in \mathbb{Z}^2 $. Then, up to a subsequence, $ \yn \equiv y_0 $. Let $ y_0 \neq 0 $ and consider $ r $ the line passing through the origin and $ y_0 $. Then, define $ \Omega_0 $ as the open connected region between $ r $ and one of the axis, such that the angle between $ r $ and the axis is $ \leq \frac{\pi}{2} $ and $ \Omega = \Omega_0 \cap B_\delta $.  Take $ \delta > 0 $ such that $ \delta < |y_0| $. So, for $ x\in \Omega $, we have that $ |x-y_0| > |y_0| $. Therefore, by the Mean Value Theorem, there exists $ x_\delta \in \Omega $, satisfying
\begin{align*}
||\un||_{\ast}^{2} & = \intR \ln(1+|x-\yn|)\until^2 (x) dx \\
& \geq |\Omega| \ln(1+|x_\delta -\yn|)\until^2 (x_\delta ) \\
& = C_1 \ln(1+|x_\delta -y_0|) \geq C_1 \ln(1+|y_0|) = C_1 \ln(1+|\yn|) ,
\end{align*}
for $ C_1 > 0 $. If $ y_0 =0 $, $ \until = \un $ and the result follows immediately from (\ref{eq28}). So, in any of the cases, there exists $ C_1 >0 $ such that
$$
||\un||_{\ast}^{2}= \intR \ln(1+|x -\yn|)\until^2 (x) dx \geq C_1 \ln(1+ |\yn|) \ , \forall \ n \in \N .
$$
Now, from (\ref{2.3}), we have
$$
||\until||_{\ast}^{2} = \intR \ln (1+|x+\yn|) \un^2(x) dx \leq ||\un||_{\ast}^{2} + \ln (1+|\yn|) ||\un||_{2}^{2} .
$$
From this and (\ref{2.3}), since every norm is weakly lower semicontinuos, $ \un \RA u $ in $ L^2(\Rdois) $ and $ \until \rightharpoonup u $ in $ X $, follows that
\begin{align*}
||(-\yn) \ast u||_{\ast}^{2}&=\intR \ln(1+|x -\yn|)u^2 (x) dx \\
& \leq ||u||_{\ast}^{2}+\ln(1+|\yn|)||u||_{2}^{2} \\
& = ||\un||^{2} + ||\un||_{\ast}^{2} + 2 \ln (1+|\yn|) ||\un||_{2}^{2} \\
& = ||\un||^{2}  + ||\un||_{\ast}^{2} (1 + C_2||\un||_{2}^{2}) \\
& \leq ||\un||^{2}  + C_3 ||\un||_{\ast}^{2} \leq C_4 ||\un||_{X}^{2}
\end{align*}
for $ n\in \N $ and $ C_4 >0 $. Consequently, there exists a constant $ C_5 >0 $ such that, after passing to a subsequence, we have, for all $ n\in \N $,
\beq\label{eq29}
||(-\yn) \ast u||_X^{2} =||u||^2 + ||(-\yn) \ast u||_{\ast}^{2} \leq ||\un||^2 +C_4 ||\un||_{X}^{2} \leq C_5 ||\un||_{X}^{2} .
\eeq

From this and (\ref{eq28}),
$$
|I'(\until)(\until - u)|\leq (1+ \sqrt{C_5})||I'(\un)||_{X'}||\un||_X \RA 0 ,
$$
as $ n\RA + \infty $. Hence, Claim 2 is proved. \\

\textbf{Claim 3:} $ \ds_{\Rdois} f(\until)(\until - u) dx \RA 0 $, as $ n\RA + \infty $. \\

Since $ ||\until||= ||\un|| $, we know that $\ds_{\Rdois} \left(e^{r\alpha \until^2}-1\right) dx \leq C_\alpha $, for $ r\sim 1 $. Thus, by (\ref{eq2}), Hölder and generalized Hölder Inequalities, and the facts that $ q>4 $ and $ \until \RA u $ in $ \Ls $, for $ s \geq 2 $, follows
\begin{align*}
\left| \intR f(\until)(\until - u) dx \right| 
& \leq ||\until||_2 ||\until - u||_{2} + b_1 C_{\alpha}^{\frac{1}{r}}  ||\until||_{r'(q-1)}^{q-1}||\until - u||_{r''} \RA 0  ,
\end{align*}
where $ \frac{1}{r'}+\frac{1}{r''}+\frac{1}{r}=1 $, with $ r\sim 1 $ and $ r'' \geq 2 $, proving claim 3. From this claim, 
\begin{align*}
o(1) & = I'(\until)(\until - u ) \\
& = o(1) + ||\until ||^{2} - ||u||^2 + \dfrac{1}{4}V_{0}'(\until)(\until - u) - \intR f(\until) (\until - u ) dx \\
& = o(1) + ||\until ||^{2} - ||u||^2 + \dfrac{1}{4} [ V_{1}'(\until)(\until - u) - V_{2}'(\until)(\until - u)] .
\end{align*}
Observe that
$$
\left|V_{2}'(\until)(\until - u)\right| = |B_2 (\until^2 , \until (\until - u))| \leq ||\until||_{\frac{8}{3}}^{3}||\until - u ||_{\frac{8}{3}} \RA 0
$$
and
$$
\dfrac{1}{4}V_{1}'(\until)(\until - u)=B_1 (\until^2 , \until (\until - u)) = B_1 (\until^2 ,  (\until - u)^2) + B_1(\until^2 , u (\until - u)) .
$$
Also, since $ (\until) $ is bounded in $ X $, from Lemma \ref{l212}, $ B_1(\until^2 , u (\until - u))\RA 0 $. So, from the above estimations 
$$
o(1) =||\until||^{2}-||u||^2 +  B_1 (\until^2 ,  (\until - u)^2) + o(1) \geq ||\until||^{2}-||u||^2  + o(1) .
$$
Hence, $ ||\until||\RA ||u|| $ and $B_1 (\until^2 ,  (\until - u)^2)\RA 0$. Using the first one and the fact that $ \until \rightharpoonup u $ in $ \Hum $, we obtain that $ ||\until - u||\RA 0 $ and using the second one and the Lemma \ref{l21}, we have $ ||\until - u ||_{\ast}\RA 0 $. Therefore, $ \until \RA u $ in $ X $. 

Finally, we need to show that $ u $ is a critical point for $I $. Let $ v\in X $. Repeating the above arguments used to obtain (\ref{eq29}), we can see that
$$
||(-\yn)\ast v||_X \leq C_5 ||\un||_{X} \ , \forall \ n \in \N , C_5 >0 .
$$
So, by (\ref{eq15}),
$$
|I'(u)v|=\lim\limits_{n\RA + \infty}|I'(\until) v |\leq C_6  \lim\limits_{n \RA + \infty} ||I'(\un)||_{X'} ||\un||_{X} =0 ,
$$
proving the proposition. 

\cqd

\bo \label{obs9}
Note that, if we substitute the hypothesis $ I(\un)\RA c_{mp} $ by $ I(\un)\leq c_{mp} $, for all $ n\in \N $, then we can prove Lemma \ref{l29} and Propositions \ref{p21} and \ref{p31} by following the same steps. It happens because they essentially need the boundedness  from above of the sequence $(I(\un))\subset \R$ by the value $c_{mp}$. 
\eo

\noindent \textbf{Proof of Theorem (\ref{t11})}: From Lemma \ref{l29}, equations (\ref{eq8}) and (\ref{eq15}) and Proposition \ref{p31}, we conclude that there exists a sequence $(\until)\subset X$ such that $ \until \RA u_0 $ in $ X $, $ I(u_0)=c_{mp} $ and $ u_0 $ is a nontrivial critical point for $ I $ in $ X $. Proving (i).

In the following, define $ K=\{u\in X\setminus\{0\} \ ; \ I'(u)=0 \}$. Note that, from (i), $ u_0\in K $ and thus $ K\neq \emptyset $. Hence, we can consider $ (\un)\subset K $ such that
$$
I(\un)\RA c_g = \inf\{ I(v) \ ; \ v\in X\setminus \{0\} , I'(v)=0\} = \inf\limits_{v\in K}I(v) \in [-\infty , c_{mp}] .
$$
Note that, as $ (\un)\subset K $, by definition of $ K $, $ I'(\un)=0 $, for all $ n\in \N $. Also, from Lemma \ref{l26}, $J(\un)=\frac{5}{2}I'(\un)(\un) - P(\un) = 0$, for $ n\in \N $. Next, lets split the proof into two cases. 

\noindent \textbf{Case 1:} If $ c_g = c_{mp} $,  then $ (\un) $ satisfies (\ref{eq15}) and we already know that there exists a nontrivial critical point $ u\in X $ for $ I $ with $ I(u)=c_{mp}=c_g $. 

\noindent \textbf{Case 2:} If $ c_g < c_{mp} $, once $ I(\un)\RA c_g $, passing to a subsequence if necessary, we can assume $ I(\un)\leq c_{mp} $, for all $ n\in \N $. So, it is obvious that henceforth, we are under Remark \ref{obs9} assertion. 

From the definition of $ K $ and (\ref{eq8}), we infer $ \liminf\limits_{n} ||\un||\geq \rho > 0 $. Consequently, by Proposition \ref{p31}, there exists $ (x_n)\subset \Rdois $ and $ u\in X\setminus \{0\} $ a critical point to $ I $ such that $ \until = x_n \ast \un \RA u  $ in $ X $. 

Therefore, $ u\in K $ and 
$$
I(u)=\lim\limits_{n\RA + \infty}I(\until) = \lim\limits_{n\RA + \infty}I(\un) = c_g .
$$
In particular, we conclude that $ c_g > -\infty $. 
\cqd

\section{Proof of Theorem 1.2}

The goal of this section is to prove Theorem \ref{t12}. As some results are simple adaptations of the previous one on sections 2 and 3, we will only sketch their proofs highlighting the differences. Henceforth, we fix a closed subgroup $ G $ of $ O(2) $ and we let $ \tau : G \RA \{-1, 1 \} $ be a group homomorphism. We also consider the action of $ G $ on $ X $, $ \ast $, and the invariant subspace $ X_G \neq \{0\}$, defined previously. 

As in the previous part, we seek for critical points of $ I $ on $ X_G $. Under the principle of symmetric criticality (\cite{[23]}, Theorem 1.28)  , we know that any critical point of the restriction of $ I $ to $ X_G $, which we will also denote by $ I $, is a critical point of $ I $. 

Once again, it is easy to see that $ 0<b_0 \leq c_{mp, G} < \infty $.

\bl\label{l42}
There exists a sequence $ (\un)\subset X_G $ such that
\beq\label{eq31}
I(\un)\RA c_{mp, G} \ , \ ||I'(\un)||_{X_{G}'}(1+||\un||_{X})\RA 0 \textrm{ \ \ and \ \ } J(\un)\RA 0 \ , \textrm{ \ as \ } n\RA + \infty .
\eeq
\el
\textbf{Proof:} Here we follow the same steps as Lemma \ref{l29} and start defining the Banach space $ \tilde{X}_G = \R \times X_G $ endowed with the norm $ ||(s, v)||_{\tilde{X}_G}=(|s|^2 + ||v||_{X}^{2})^{\frac{1}{2}} $, for $ s\in \R $ and $ x\in \Rdois $. 

Again, we consider the function $ \varphi : \tilde{X}_G \RA X_G $ defined by $ \varphi(s, v)(x)=e^{\frac{5}{2}s}v(e^s x) $, for $ s\in \R $, $ x\in \Rdois $ and $ v\in X_G $. Here it is important to verify that $ \varphi $ is well-defined. In fact, given $ v\in X_G $
\begin{align*}
[A\ast \varphi(s, v)](x) & =\tau(A)\varphi(s, v)(A^{-1}x) =\tau(A)e^{\frac{5}{2}s}v(e^s A^{-1}x) \\
&=e^{\frac{5}{2}s}\tau(A)v(A^{-1}(e^{s}x))=e^{\frac{5}{2}s}[A\ast v](e^s x)\\
&= e^{\frac{5}{2}s}v(e^s x) = \varphi(s, v) .
\end{align*}
Therefore, $ \varphi(s, v)\in X_G $ for $ s\in \R $ and $ v\in X_G $. To finish the proof we just need to repeat the steps of Lemma \ref{l29} proof, so we omit it here. 

\cqd

\bo\label{obs41}
Note that, for $ A\in G\subset O(2) $, from $ |\det A| =1 $, $ \tau(A)= \pm 1 $ and $ (f_1 ') $, we have that $ I $ is invariant under the action of $ G $ and $ I'(v)(w)=0 $, for $ v\in X_G $ and $ w\in X_{G}^{\perp} $. Hence,
$$
||I'(v)||_{X'}=\sup\limits_{w\in X\setminus \{0\}}|I'(v)(w)| = \sup\limits_{w\in X_G \setminus \{0\}}|I'(v)(w)|=||I'(v)||_{X_{G}'} , 
$$
for all $ v\in X_G $, since $ X=X_{G}\oplus X_{G}^{\perp} $. Consequently, we can rewrite (\ref{eq31}) as
\beq\label{eq32}
I(\un)\RA c_{mp, G} \ , \ ||I'(\un)||_{X'}(1+||\un||_{X})\RA 0 \textrm{ \ \ and \ \ } J(\un)\RA 0 \ , \textrm{ \ as \ } n\RA + \infty .
\eeq
\eo

\bl\label{l43}
Let $ (\un)\subset X_G $ satisfying (\ref{eq32}). Then, $ (\un) $ is bounded in $ \Hum $ and for a sufficiently small $ \rho_0 >0 $, we have $ \limsup\limits_{n}||\un|| < \rho_0 $. 
\el

The next result is the key to prove Theorem \ref{t12} and it has a quite different proof than its section 3 analogues, so we will present its full version. To do that, we need to define the set $\textrm{ \ Fix }(G)=\{x\in \Rdois \ ; \ Ax=x , \forall \ A \in G\}\subset \Rdois $.

\bp\label{p41}
Let $ q>4 $ and $ (\un)\subset X_G $ satisfying (\ref{eq32}). Then, after passing to a subsequence, one between the two alternatives happens:

\textbf{(a)} $ ||\un||\RA 0 $ and $ I(\un)\RA 0 $, as $ n\RA + \infty $. 

\textbf{(b)} There exists points $ y_n \in \textrm{ \ Fix }(G) $ such that $ y_n \ast \un \RA u $ in $ X $, for some critical point $ u\in X_G \setminus \{0\} $ of $ I $.
\ep
\textbf{Proof:} Suppose that (a) does not occur. Then, by Proposition \ref{p31}, passing to a subsequence if necessary, $ \tilde{y}_n\ast \un \RA \util $ in $ X $, for suitably $ \tilde{y}_n\in \Rdois $ and a nontrivial critical point $ \util \in X $ of $ I $. \\

\noindent \textbf{Claim:} $ \omega = \sup\{|A\tilde{y}_n - \tilde{y}_n| \ ; \ A\in G , n\in \N\} < \infty $. \\

\noindent Set $ z_n = -\tilde{y}_n $, for $ n\in \N $. Let $ (A_n)\subset G $ be any sequence. As, for $ n \in \N  $, $ \un \in X_G $, $ \ast $ is linear and $ ||\cdot ||_2 $ is invariant under the action of $ G $, we have that
\begin{align}
||A_n \ast (z_n \ast \util) - z_n \ast \util||_{2} 
& \leq ||A_n \ast (z_n \ast \util - \un)||_{2} + ||\un - z_n \ast \util||_{2} \nonumber \\
& = 2||z_n \ast \util - \un ||_{2}= 2||\util - \tilde{y}_n \ast \un||_{2}\RA 0 , \label{eq33}
\end{align}
as $ n\RA + \infty $. Now, set $ \vn = A_n \ast \util $ and $\xi_n = A_n z_n -z_n$, for $ n\in \N $. Then, considering $ x-z_n = y $, one can easily see that $ (\xi_n \ast \vn)(y) = \tau(A_n) \util ( A_{n}^{-1}(y-\xi_n))$ holds and, by changing of variables,  
\begin{align}
||A_n \ast (z_n \ast \util) - z_n \ast \util||_{2}^{2} - 2||\util||_{2}^{2} &= -2\tau(A_n)\intR \util(A_{n}^{-1}x-z_n)\util(x-z_n) dx \nonumber \\
& = -2\tau(A_n)\intR \util(A_{n}^{-1}(y-\xi_n))\util(y) dy \nonumber \\
& = -2 \intR (\xi_n \ast \vn) \util dx . \label{eq34}
\end{align}
Since $ G $ is a compact set, as a closed subgroup of $ O(2) $, passing to a subsequence if necessary, $ A_n \RA A \in G $. Then, setting $ v=A\ast \util $, $ \vn = A_n \ast \util \RA v $ in $ X $ and, consequently, 
\beq\label{eq35}
||\xi_n \ast (\vn -v)||_{2}=||\vn - v ||_{2}\RA 0 , \textrm{ \ as \ } n \RA + \infty .
\eeq
Thus, if we put $ \beta_n = ||A_n \ast (z_n \ast \util) - z_n \ast \util||_{2}\RA 0 $, from equations (\ref{eq32})-(\ref{eq35}) and Hölder inequality, we deduce that
\begin{align*}
\left|\intR (\xi_n \ast v) \util dx - \intR |u|^2 dx \right| & = \left| \intR (\xi_n \ast v) \util dx - \dfrac{1}{2}\beta_n - \intR (\xi_n \ast \vn) \util dx\right| \\
& \leq ||\xi_n \ast (\vn - v)||_2 ||\util||_2 + \dfrac{1}{2}\beta_n\RA 0 ,
\end{align*}
as $ n\RA + \infty $. Therefore, 
\beq\label{eq36}
\lim\limits_{n\RA + \infty} \intR (\xi_n \ast v) \util dz =||u||_{2}^{2}>0 .
\eeq
Based on this, we conclude that $ (\xi_n) $ remains bounded, since, otherwise, after passing to a subsequence, $ \xi_n \ast v \rightharpoonup 0$ in $ \Ldois $. Finally, $ |A_n \tilde{y}_n - \tilde{y}_n|=|A_n z_n - z_n| =|\xi_n| $, for $ n \in \N $, and the claim follows. 

Lets substitute $ \tilde{y}_n $ by $ \yn = \dfrac{1}{\mu(G)}\ds_{G}A\tilde{y}_n d\mu(A)\in \textrm{ Fix} (G) $, for $ n\in \N $, where $ \mu $ is the Haar measure of $ G $ and $ \tilde{y}_n $ is fixed, for each $ n\in \N $. From Claim 1, 
$$
|y_n - \tilde{y}_n| \leq \left| \dfrac{1}{\mu(G)}\ds_{G}A\tilde{y}_n d\mu(A) - \tilde{y}_n \right| \leq \dfrac{1}{\mu(G)}\ds_{G}|A\tilde{y}_n - \tilde{y}_n| d\mu(A) \leq \omega .
$$
Thus, $(y_n -\tilde{y}_n)$ is bounded in $ \Rdois $ and, passing to a subsequence if necessary, $ y_n -\tilde{y}_n\RA r\in \Rdois $. Consequently, setting $ u=r\ast \util $, $ y_n \ast \un \RA u  $ in $ X $. 
Indeed, 
it follows noticing that
\begin{eqnarray*}
\|y_n \ast u_n -r\ast \tilde{u}\|_{X} &\leq& \|(y_n  -(y_n-\tilde{y}_n))\ast u_n- \tilde{u}\|_{X}+\| (y_n-\tilde{y}_n)\ast \tilde{u}\ -r\ast \tilde{u}\|_{X}+o(1)\\
&=& o(1),\ \mbox{as}\ n \rightarrow \infty.
\end{eqnarray*}
To finish the proof we need to verify that such function $ u $ is on $ X_G $. In fact, for $ A\in G $ and $ x\in \Rdois $, we have, for $ y_n\in \textrm{ Fix} (G)  $,
\begin{align*}
[A\ast u](x) & = \lim\limits_{n\RA + \infty} A\ast (y_n \ast \un) 
 = \tau(A)\lim\limits_{n\RA + \infty}\un(A^{-1}(x-y_n)) \\
& = \lim\limits_{n\RA + \infty}[y_n \ast (A \ast \un)](x) = \lim\limits_{n\RA + \infty} [y_n \ast \un ](x) = u(x) .
\end{align*}
Therefore, $ u\in X_G $, proving that (b) occurs. 
\cqd

\bo\label{obs42}
Once again, we can substitute the hyphotesis that $I(\un)\RA c_{mp, G}$ by \linebreak $ I(\un)\leq c_{mp, G} $, for all $ n\in \N $, and prove the same above statements.
\eo

\textbf{Proof of Theorem \ref{t12}:} From Proposition \ref{p41} and the fact that $ c_{mp, G}>0 $, there exists a critical point $ u_0 \in X_{G}\setminus \{0\} $ of $ I $ with $ I(u_0)= c_{mp, G}$, proving item (i). \\

In particular, $ K=\{ u \in X_{G}\setminus \{0\} \ ; \ I'(u)=0\} \neq \emptyset $, since $ u_0\in K $. Consider a sequence $(\un)\subset K $ such that 
$$
I(\un)\RA c_{g, G}=\inf\limits_{v\in K}I(v) \in [-\infty , c_{mp, G}] .
$$
From the definition of $ K $ and the convergence, arguing as in the proof of Theorem \ref{t11}, we have that 
$$
I(\un)\leq c_{mp, G} \ , \ ||I'(\un)||_{X'}(1+||\un||_{X})\RA 0 \textrm{ \ \ and \ \ } J(\un)\RA 0 \ , \textrm{ \ as \ } n\RA + \infty ,
$$
so, in view of Remark \ref{obs42}, we can apply the results of this section to sequence $ (\un) $. Also, from (\ref{eq8}), $ ||\un||\geq \rho > 0 $, for all $ n\in \N $, since otherwise $ n_\beta =0 $. Thus, $ \liminf\limits_{n}||\un|| \geq \rho > 0 $. 

Consequently, by Proposition \ref{p41}, there exists points $ y_n \in \textrm{ \ Fix}(G) $ and a nontrivial critical point $ u\in X_G \setminus\{0\} $ of $ I $ such that, passing to a subsequence, $ y_n \ast \un \RA u $ in $ X $. Hence, $ I'(u)=0 $, which implies that $ u\in K $ and
$$
I(u)=\lim\limits_{n\RA + \infty} I(y_n \ast \un) =\lim\limits_{n\RA + \infty} I(\un) = c_{g, G} .
$$
Therefore, item (ii) occurs and, in particular, $ c_{g, G} > -\infty $.

\cqd



\noindent \textbf{Acknowledgements:} This work was perfomed and completed during the  first author PhD graduate course at Federal University of São Carlos. \\

\end{document}